\documentclass[10pt]{article}
\usepackage{amsmath,amssymb,pictex,placeins}
\usepackage[margin=20pt,font=small,labelfont=bf]{caption}
\usepackage{afterpage}
\newcount\size
\newcount\msize
\newbox\bottomedge
\newbox\topedge
\newbox\interiorone
\newbox\interiortwo
\newbox\interiorthree
\newdimen\gridsize
\def\begincolor#1{}

\begin{document}

\title{\bfseries The Domination Number of
  $C_n\square P_m$ for $n\equiv 2\pmod{5}$}
\author{%
    David R. Guichard\\%
    Whitman College\\%
    345 Boyer Ave\\%
    Walla Walla, WA 99362\\%
}
\date{}
\maketitle
\thispagestyle{empty}

\begin{abstract} 
We use a dynamic programming algorithm to establish a new lower bound
on the domination number of complete cylindrical grid graphs of the
form $C_n\square P_m$, that is, the Cartesian product of a path and a
cycle, when $n\equiv 2\pmod{5}$, and we establish a new upper bound
equal to the lower bound, thus computing the exact domination number
for these graphs.

%
\end{abstract}

\section{Introduction} A set $S$ of vertices in a graph $G=(V,E)$ is
called a {\it dominating set\/} if every vertex $v\in V$ is either in
$S$ or adjacent to a vertex in $S$.  The domination number of $G$,
$\gamma(G)$, is the minimum size of a dominating set.

Let $P_m$ denote the path on $m$ vertices and $C_n$ the cycle on $n$
vertices; the {\it complete cylindrical grid graph\/} or
{\it cylinder} is the product
$C_n\square P_m$. That is, if we denote the vertices of $C_n$ by
$u_1,u_2,\ldots,u_n$ and the vertices of $P_m$ by $w_1,\ldots,w_m$,
then $C_n\square P_m$ is the graph with vertices $v_{i,j}$, $1\le
i\le n$, $1\le j\le m$, and $v_{i,j}$ adjacent to $v_{k,l}$ if
$i=k$ and $w_j$ is adjacent to $w_l$ or if
$j=l$ and $u_i$ is adjacent to $u_k$.
It will  be useful to think of this graph as $P_n\square P_m$,
with the edge paths of length $m$ glued together, that is, connected
with new edges. We assume throughout that $m\ge 2$ and $n\ge 3$.

P. Pavli{\v c} and J. {\v Z}erovnik
\cite{pavlic-zerovnik:dom-no-cyl-graphs} established upper bounds for
the domination number of $C_n\square P_m$, and 
Jos\'e Juan Carre{\~n}o et al 
\cite{carreno-et-al:lower-bound-dom-no-cyl-graphs}
established non-trivial lower bounds. For $n\equiv 0\pmod 5$ the
bounds agree, so the domination number is known exactly.
In \cite{guichard:domination_bounds_cylinder},
we
improved the lower bounds, except of course in the case that
$n\equiv 0\pmod 5$. Using the same techniques, we improve the lower
bound slightly when $n\equiv 2\pmod 5$, and we show that this value is
the true domination number by also slightly improving the known upper
bound. 

\section{Lower bound} 
We summarize the technique for computing a lower bound; details are in
\cite{guichard:domination_bounds_cylinder}.
A vertex in $C_n\square P_m$ dominates at most five vertices,
including itself, so certainly $\gamma(C_n\square P_m)\ge nm/5$. If we could
keep the sets dominated by individual vertices from overlapping, we
could get a dominating set with approximately $nm/5$ vertices, and
indeed we can arrange this for much of the graph, with the exception
of the ``edges'', that is, the two copies of $C_n$ in which the
vertices have only 3 neighbors, and except when $n\equiv 0\pmod 5$, we
run into some trouble where the edge paths of $P_n\square P_m$ are
connected to form $C_n\square P_m$.

Suppose $S$ is a subset of the vertices of $C_n\square P_m$. Let $N[S]$ be
the set of vertices that are either in $S$ or adjacent to a member of
$S$, that is, the vertices dominated by $S$.  Define the {\it wasted
domination\/} of $S$ as $w(S)=5|S|-|N[S]|$, that is, the number of
vertices we could dominate with $|S|$ vertices in the best case, less
the number actually dominated. When $S$ is a dominating set,
$|N[S]|=mn$, and if $w(S)\ge L$ then $|S|\ge (L+mn)/5$. Our goal now
is to find a lower bound $L$ for $w(S)$.

\begin{figure}[htb]
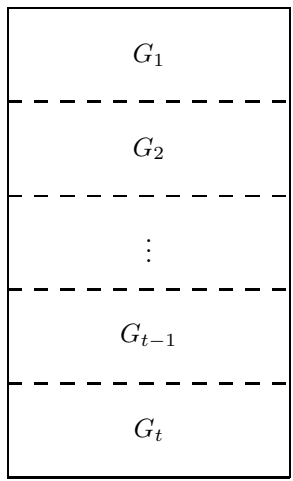

\hbox to \hsize{\hss
\beginpicture
\setcoordinatesystem units <2.5mm,2.5mm>
\setplotarea x from 0 to 15, y from 0 to 25
\axis left /
\axis right /
\axis top /
\axis bottom /
\setdashes
\putrule from 0 5 to 15 5
\putrule from 0 10 to 15 10
\putrule from 0 15 to 15 15
\putrule from 0 20 to 15 20
\putrule from 0 25 to 15 25
\put {$G_1$} at 7.5 22.5
\put {$G_2$} at 7.5 17.5
\put {$\vdots$} at 7.5 12.5
\put {$G_{t-1}$} at 7.5 7.5
\put {$G_t$} at 7.5 2.5
\endpicture\hss}
\caption{Partitioned cylinder.}
\label{fig:partitioned}
\end{figure}

Suppose a cylinder $C_n\square P_m$ is partitioned into
subgraphs as indicated in Figure~\ref{fig:partitioned}, where each $G_i$
is a subgraph $C_n\square P_{m_i}$. (We will refer to $G_1$ and $G_t$
as edge subgraphs, and the others as interior subgraphs.) 
Let $S$ be a
dominating set for $G$ and $S_i=S\cap V(G_i)$.  Then
\begin{equation}\label{eq:fundamental inequality}
w(S)\ge \sum_{i=1}^t w(S_i).
\end{equation}
Note that in computing $w(S_i)$ we consider $S_i$ to be a subset of $V(G)$,
not of $V(G_i)$ (this affects the computation of $N[S_i]$). To verify
the inequality, note that the following inequalities are equivalent:
\begin{align*}
  w(S)&\ge \sum_{i=1}^t w(S_i)\\
  5|S|-|N[S]|&\ge \sum_{i=1}^t (5|S_i|-|N[S_i]|)\\
  5|S|-|N[S]|&\ge \sum_{i=1}^t 5|S_i|-\sum_{i=1}^t |N[S_i]|\\
  |N[S]| &\le \sum_{i=1}^t |N[S_i]|.  
\end{align*}
The last inequality is satisfied, since each vertex in $N[S]$ is
counted at least once by the expression on the right. 

Note that $S_i$ is a set that dominates all the vertices of $G_i$
except possibly some vertices in the top or bottom row of $G_i$ (or in
the cases of $G_1$ and $G_t$, in the bottom row and top row,
respectively).  Let us say that a set that dominates a
cylinder $G$, with the exception of some vertices on the
top or bottom edges, {\it almost dominates\/} $G$.  Given a 
cylinder $H=C_n\square P_{m_i}$ (namely, one of the $G_i$),
What we want to know is the value of
\begin{equation}\label{eq:desired minimum}
  \min_A w(A),
\end{equation}
taking the minimum over sets $A$ that almost dominate
$H$ and computing $w(A)$ as if $A$ were a subset of a larger graph
$C_n\square P_{m_i+2}$ in which $H$ occupies the middle $m_i$ rows, or in
the case of $G_1$ or $G_t$, $A$ is a subset of $C_n\square P_{m_i+1}$ in
which $H$ occupies the top $m_i$ rows.  If we can compute this minimum
for (small) fixed $m_i$ and any $n$, we can choose $G_1$ through $G_t$
with a small number of rows and get lower bounds on $w(S_i)$ for any
dominating set $S$ of the original $C_n\square P_m$.

\section{The algorithm}

We use a dynamic programming algorithm to compute these lower bounds
on $w(S_i)$. This is done by computing smallest almost-dominating sets
for $C_i\square P_j$, $i=1,2,3\ldots$, until periodicity is detected.
That is, we test for a point after which the minimum wasted domination
of $C_i\square P_j$ is that of $C_{i-p}\square P_i$ plus a constant
$q$.  (Livingston and Stout~\cite{ls:constant-time-dominating-sets}
and Fisher~\cite{dcf:domination-grid-graphs} independently thought of
looking for this sort of periodicity.)  In
\cite{guichard:domination_bounds_cylinder}, we found that the minimum
wasted domination for $G_1=C_n\square P_{10}$ is $n$. For an interior
$G_i=C_n\square P{10}$, the minimum wasted domination is 0, 6, 5, 9,
or 6 as $n$ is 0, 1, 2, 3, or 4 $\pmod 5$.  For $n\equiv 2\pmod{5}$,
this means that for a dominating set $S$ of $C_n\square P{m}$, we have
$$5|S|-mn = w(S) \ge 2n+5\lfloor\frac{m-20}{10}\rfloor,$$
so that 
$$|S|\ge \frac{(m+2)n}{5}+\lfloor \frac{m-20}{10}\rfloor.$$
Since $m$ need not be a multiple of 10, this calculation effectively
ignores one of the $G_i$, namely, one in which the number of rows is
less than 10. When the number of rows is too small (less than 5, as it
turns out), the resulting lower bound is too small to be useful,
namely, 0. So we computed the corresponding minimum wasted domination
for an interior $G_i=C_n\square P_j$, for $5\le j\le 9$. We find that
the minima are 2, 3, 3, 4, and 4, respectively. This means, for example, that if
$m\equiv 8\pmod{10}$,
\begin{align*}
  5|S|-mn &\ge 2n+5\lfloor\frac{m-20}{10}\rfloor+4\\
  |S|&\ge\frac{(m+2)n}{5}+\lfloor \frac{m-20}{10}\rfloor+\frac{4}{5}.\\
\end{align*}
To handle $m\equiv j\pmod{10}$, $1\le j\le 4$, we eliminate one
interior $G_i$ with 10 rows, and insert two new subgraphs, with
5 and 6 rows, 6 and 6 rows, 6 and 7 rows, or 6 and 8 rows,
respectively. For example, for $m\equiv 2\pmod{10}$, we get
$$|S|\ge \frac{(m+2)n}{5}+\lfloor \frac{m-20}{10}\rfloor-1+\frac{3}{5}+\frac{3}{5}
$$

When we complete this task, the resulting lower bounds still turn out
to be a bit too low for some values of $m$. To improve the bounds, we
compute the minimum wasted domination for edge subgraphs $G_1$ and
$G_t$ with 11 rows, instead of 10. We find that the minimum wasted
domination for $C_n\square P_{11}$ is $n+2$, when $n\equiv
2\pmod5$. The resulting lower bound on the domination number of
$C_n\square P_m$ becomes

$$|S|\ge \frac{(m+2)n+4}{5}+\lfloor \frac{m-22}{10}\rfloor+c,
$$
where $c$ is the necessary addition as described above when $m-22$ is
not divisible by 10. For example, when $m\equiv 0\pmod{10}$, we get
$$|S|\ge \frac{(m+2)n+4}{5}+\lfloor \frac{m-22}{10}\rfloor-1+2\cdot\frac{3}{5},
$$
since now $m-22\equiv 8\pmod{10}$.

The resulting lower bounds for $m=10i+j$ and $n=5k+2$ are:
\begin{align}
  \begin{split}
  \label{eq:lb}
  j=0\colon\;&10ik+5i+2k\\
  j=1\colon\;&10ki+5i+3k\\
  j=2\colon\;&10ik+5i+4k+1\\
  j=3\colon\;&10ik+5i+5k+1\\
  j=4\colon\;&10ik+5i+6k+2\\
  j=5\colon\;&10ik+5i+7k+2\\
  j=6\colon\;&10ik+5i+8k+3\\
  j=7\colon\;&10ik+5i+9k+3\\
  j=8\colon\;&10ik+5i+10k+4\\
  j=9\colon\;&10ik+5i+11k+4\\
  \end{split}
\end{align}
With a little bit of algebraic manipulation, we find that all of these
bounds can be expressed as
$$
\left\lceil\frac{\frac{4n+2}{10}m+\frac{4n-13}{5}}{2}\right\rceil.
$$
These bounds apply when $n\ge 32$, $n\equiv 2\pmod{5}$; and when $m\ge
27$. The first restriction is due to when the periodicity described
earlier first appears; the second because the bounds are based on two
edge subgraphs with 11 rows and at least one interior subgraph with at
least 5 rows.

Crevals~\cite{crevals:domination_cylinder_graphs} computes exact
values for $m\le 22$ and all $n\ge m$, and also for
$n\le 30$ and all $m\ge n$. Our bounds agree with Crevals' exact values
for $16\le m\le 22$, and for $n\le 30$ (remember that our bounds are for
$n\equiv 2\pmod5$ only). Note well that our bounds agree with the
Crevals values when $n<32$ and when $16\le m<22$, even though our derivation
of the lower bounds does not apply to these values.

This gave us some confidence that our lower bounds would prove correct
in all cases except when $m\le 15$. To show this, we will need to find
upper bounds that match our lower bounds, and we also need to fill in
the gap $22<m<27$ not covered by either our lower bounds or the
Crevals values. So the first order of business is to provide lower
bounds for these missing values.

This we can do in a way similar to what we have already done, but we
need to compute minima for edge subgraphs of the form $C_n\square
P_{12}$ and $C_n\square P_{13}$. When we do this, we find that the
minimum wasted domination for both is $n+3$, when $n\equiv
2\pmod5$. Then we compute lower bounds as before, using two edge
subgraphs and no interior subgraphs. To get bounds for $m=23,24,25,26$
we use subgraphs with 11 and 12 rows; 12 and 12 rows; 12 and 13 rows;
and 13 and 13 rows; respectively.  We find that the resulting lower
bounds are given by the formulas we already have. For example, for
$m=25$ and $n=5k+2$ we have
$$|S|\ge \frac{25n}{5}+\frac{2n+6}{5}=27k+12,$$
which agrees with the $j=5$ value in equation~\ref{eq:lb},
$10ik+5i+7k+2$, since $i=2$.

\section{Upper bound}

The best known upper bound, due to Pavli{\v c} and {\v Z}erovnik
\cite{pavlic-zerovnik:dom-no-cyl-graphs}, for $n\equiv 2\pmod5$ is
$\frac{(m+2)n}{5}+\frac{1}{10}(m+2)$. This is slightly larger than our
lower bound. We need to show that $C_n\square P_m$ can in fact be
dominated by the number of vertices given by our lower bound. We can
easily modify our computer programs to provide almost-dominating sets
of the component subgraphs $G_i$ with the minimum possible wasted
domination. We would then hope that these subgraphs can be pieced
together in a way that dominates the entire cylinder; we can attempt
to do this for a particular value of $n$, say $n=32$. To then extend
the upper bound to all values of $n$, we would need these sample
graphs to exhibit a periodicity that allows them to be extended to any
$n$. With a little bit of trial and error, we succeeded.

We discovered that the sample interior graphs that fit together nicely
are the ones with an even number of rows. So to begin, we need to show
that we can get dominating sets of cylinders with $m$ rows, for all
values of $m\bmod 10$, using only such graphs. We can potentially get
these for all even $m$ using two edge graphs of 11 rows, for $m\ge 28$
(because the graphs with fewer than 6 rows are not helpful). To get
odd values of $m\ge 27$, we use one edge graph of 11 rows and one of
10 rows. For $23\le m\le 25$, we use two edge graphs of 11, 12, or 13
rows as needed, and for $m=26$ we use two edge graphs with 10 rows and
one interior graph with 6 rows. It then suffices to provide examples
for $m=33\ldots 42$, since for $m$ larger than this it will be clear
that any number of additional interior subgraphs with 10 rows can be
added. For $m=27\ldots 32$, it will be apparent that some graphs with
$m$ between 33 and 42 can easily be modified (by removing one interior
subgraph) to verify the upper bounds for $m$ between 27 and 32.  There
was one slight surprise: for $m$ even, $m\ge 28$, one of the two edge
graphs with 11 rows must have a single vertex added to form a
dominating set.

\begin{figure}[htb]
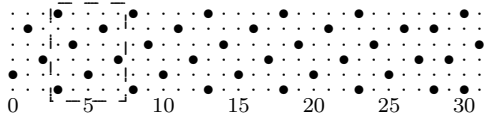

\vbox{%
  \footnotesize
  \centerline{\input graph_6x32_interior_bare_rot_0_repeat}
}
\caption{Interior subgraph $C_{32}\square P_{\the\msize}$}
\label{fig:32x6}
\end{figure}

\begin{figure}[htb]
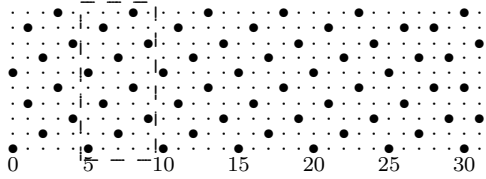

\vbox{%
  \footnotesize
  \centerline{\input graph_10x32_interior_bare_rot_0_repeat}
}
\caption{Interior subgraph $C_{32}\square P_{\the\msize}$}
\label{fig:32x10}
\end{figure}

In Figures~\ref{fig:32x6} to \ref{fig:32x10} we show the interior
graphs used to construct dominating sets. In each we show a five
column portion of the graph that can be repeated as necessary to get a
cylinder with larger $n$. We have omitted the edges in these graphs as
the structure is apparent; the small dots are vertices, the large dots
are vertices in the dominating sets.

\msize0

\msize0
\begin{figure}[htb]
\vbox{%
  \footnotesize
  \centerline{\input graph_8x32_interior_bare_rot_0_repeat}
}
\caption{Interior subgraph $C_{32}\square P_{\the\msize}$}
\label{fig:32x8}
\end{figure}

\msize0

\begin{figure}[htb]
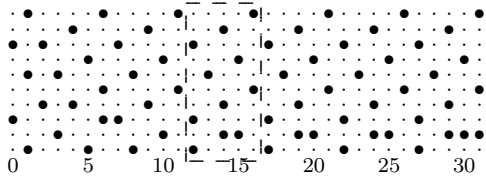

\vbox{%
  \footnotesize
  \centerline{\input graph_10x32_edge_bare_rot_0_repeat}
}
\caption{Edge subgraph $C_{32}\square P_{\the\msize}$}
\label{fig:32x10e}
\end{figure}

\begin{figure}[htb]
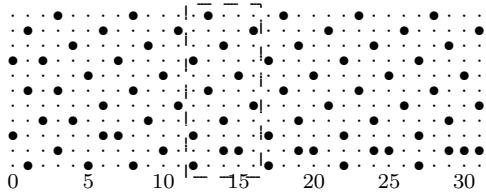

\vbox{%
  \footnotesize
  \centerline{\input graph_11x32_edge_bare_rot_0_repeat}
}
\caption{Edge subgraph $C_{32}\square P_{\the\msize}$}
\label{fig:32x11e}
\end{figure}

\begin{figure}[htb]
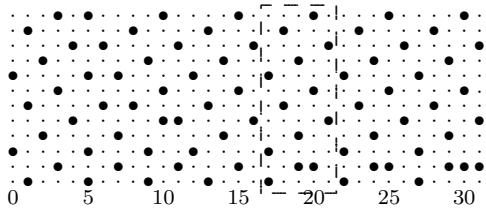

\vbox{%
  \footnotesize
  \centerline{\input graph_12x32_edge_bare_rot_0_repeat}
}
\caption{Edge subgraph $C_{32}\square P_{\the\msize}$}
\label{fig:32x12e}
\end{figure}

\begin{figure}[htb]
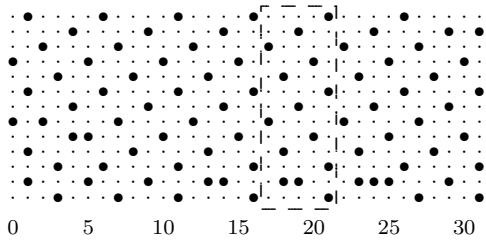

\vbox{%
  \footnotesize
  \centerline{\input graph_13x32_edge_bare_special_repeat}
}
\caption{Edge subgraph $C_{32}\square P_{\the\msize}$}
\label{fig:32x13e}
\end{figure}

In Figures~\ref{fig:32x10e} to \ref{fig:32x13e} we show the edge
graphs used to construct dominating sets. Again, in each we show a five
column portion of the graph that can be repeated as necessary to get a
cylinder with larger $n$.

\msize0

\msize0

\msize0

\msize0

Finally, we show some examples of cylinders put together from these
pieces. In the case of even $m$, as mentioned above, we need an extra
vertex in the dominating set; this vertex is shown as a star. The
extra vertex is always in the bottom row of the top edge graph.  Note
also that to make the subgraphs fit together properly, so that we get
a dominating set for the whole cylinder, some of the subgraphs must be
rotated so that they line up properly at the boundaries, which are
shown with dashed lines. (We produced in all twenty of these graphs,
covering $23\le m\le 42$.)

\begin{figure}[htb]
\vbox{%
\footnotesize\input composite_25x32_inc
}
\caption{Partitioned cylinder $C_{32}\square P_{\the\msize}$,
$|S|=\the\size$}
\label{fig:32x25}
\end{figure}

\begin{figure}[htb]
\vbox{%
\footnotesize\input composite_26x32_inc
}
\caption{Partitioned cylinder $C_{32}\square P_{\the\msize}$,
$|S|=\the\size$}
\label{fig:32x26}
\end{figure}

\begin{figure}[htb]
\vbox{%
\footnotesize\input composite_33x32_inc
}
\caption{Partitioned cylinder $C_{32}\square P_{\the\msize}$,
$|S|=\the\size$}
\label{fig:32x33}
\end{figure}

\begin{figure}[htb]
\vbox{%
\footnotesize\input composite_34x32_inc
}
\caption{Partitioned cylinder $C_{32}\square P_{\the\msize}$,
$|S|=\the\size$}
\label{fig:32x34}
\end{figure}

\begin{figure}[htb]
\vbox{%
\footnotesize\input composite_37x32_inc
}
\caption{Partitioned cylinder $C_{32}\square P_{\the\msize}$,
$|S|=\the\size$}
\label{fig:32x37}
\end{figure}

\begin{figure}[htb]
\vbox{%
\footnotesize\input composite_38x32_inc
}
\caption{Partitioned cylinder $C_{32}\square P_{\the\msize}$,
$|S|=\the\size$}
\label{fig:32x38}
\end{figure}

\FloatBarrier

To produce a dominated cylinder of arbitrary size, we first determine
which of the individual subgraphs we will need, namely, two edge
subgraphs, some number of interior subgraphs with 10 rows, and one or
two additional graphs with 6 or 8 rows. Then using the repeatable
sections in the graphs of Figures~\ref{fig:32x6} to \ref{fig:32x13e},
we expand each of these graphs to the desired number of columns, and
finally we rotate the graphs as necessary so that they fit together
properly.

We can now write down an upper bound for each $m$ and $n$. For
example, for $m\equiv 7\pmod{10}$ and $m\ge 27$, we use edge graphs
with 11 and 10 rows, an interior graph with 6 rows, and as many
additional interior subgraphs of 10 rows
as needed. Let $m=10i+7$ and $n=5k+2$. The
number of vertices contributed by the edge graph with 11 rows is
$77 + \frac{n-32}{5}12=12k+5$, since the repeatable section, from
Figure~\ref{fig:32x11e}, contributes 12 vertices to the dominating set.
Similarly, the edge graph with 10 vertices contributes
$71+\frac{n-32}{5}11=11k+5$, an interior graph with 10 rows
contributes $65+(k-6)10=10k+5$, and an interior graph with 6 rows
contributes $39+(k-6)6=6k+3$. Then the total number of vertices in a
dominating set is
\begin{align*}
(12k+5)+&(11k+5)+(6k+3)+(10k+5)\frac{m-27}{10}=\\
  &(12k+5)+(11k+5)+(6k+3)+(10k+5)(i-2)=\\
  &10ik+5i+9k+3,
\end{align*}
which matches the corresponding lower bound in equation~\ref{eq:lb}.
We repeat this process for all other $j$ in $0,\dots,9$, with
$m=10i+j$, and in each case we find that the results match the values
in equation~\ref{eq:lb}. We also do $23 \le m\le 26$ as special cases,
with no surprises.

\section{Summary}

For $m\ge 2$ and $n\ge 3$,
the domination number of $C_n\square P_m$ is now known when $n\equiv
0\pmod{5}$ (as a result of P. Pavli{\v c} and J. {\v Z}erovnik
\cite{pavlic-zerovnik:dom-no-cyl-graphs} and Jos\'e Juan Carre{\~n}o
et al \cite{carreno-et-al:lower-bound-dom-no-cyl-graphs}) and $n\equiv
2\pmod{5}$. For the latter, when $m\ge 16$ or $m\ge n$, the domination
number is
$$
\left\lceil\frac{\frac{4n+2}{10}m+\frac{4n-13}{5}}{2}\right\rceil,
$$
as shown here.
For $m<16$ and $n>m$, the domination numbers are given by the formulas
found in Crevals~\cite{crevals:domination_cylinder_graphs}.

Unfortunately, the technique used here seems unlikely to settle the
remaining cases, when $n$ is 1, 3, or $4\pmod{5}$. The lower bounds we
get are below the upper bounds of
\cite{pavlic-zerovnik:dom-no-cyl-graphs} by a small multiple of
$m$. When we look at a few sample subgraphs generated by our
algorithm, they do not fit together nicely at the boundaries (that is,
some vertices are left undominated), and so
do not give upper bounds matching our lower bounds. Indeed, it seems
quite remarkable that for $n\equiv 2\pmod5$ the subgraphs can be
combined to dominate the entire graph.

It is possible that our technique might work using somewhat larger
subgraphs, though the algorithm we used would have to be improved or
replaced with a substantially faster one. Computing the minimum
wasted domination for fewer than 13 rows was reasonably fast on an
Intel Core i7, although for 12 rows it did take a week or two. The 13
row graphs took considerably longer. We split the problem into 10
pieces by hand and ran each on a separate Intel Core i5 computer that
was otherwise idle. The
division wasn't perfect, so some finished earlier than others; the
whole process took over three months. The programs were written in C++
and compiled with the Gnu gcc compiler with optimization. 

\bibliography{guichard}

\def\noopsort#1{}
\begin{thebibliography}{1}

\bibitem{carreno-et-al:lower-bound-dom-no-cyl-graphs}
J.~J. Carre{\~n}o, J.~A. Mart{\'i}nez, and M.~L. Puertas.
\newblock A general lower bound for the domination number of cylindrical
  graphs.
\newblock {\em Bull. Malays. Math. Sci. Soc.}, 43(2):1671--1684, 2020.

\bibitem{crevals:domination_cylinder_graphs}
S.~Crevals.
\newblock Domination of cylinder graphs.
\newblock {\em Congr. Numer.}, 219:53--63, 2014.

\bibitem{dcf:domination-grid-graphs}
D.~C. Fisher.
\newblock The domination number of complete grid graphs.
\newblock manuscript.

\bibitem{guichard:domination_bounds_cylinder}
D.~R. Guichard.
\newblock A new lower bound for the domination number of complete cylindrical
  grid graphs.
\newblock {\em J. Combin. Math. Combin. Comput.}, to appear.

\bibitem{ls:constant-time-dominating-sets}
M.~Livingston and Q.~Stout.
\newblock Constant time computation of minimum dominating sets.
\newblock {\em Congr. Numer.}, 105:116--128, 1994.

\bibitem{pavlic-zerovnik:dom-no-cyl-graphs}
P.~Pavli{\v c} and J.~{\v Z}erovnik.
\newblock A note on the domination number of the cartesian products of paths
  and cycles.
\newblock {\em Kragujevac J.\ Math.}, 37:275--285, 2013.

\end{thebibliography}
\bibliographystyle{abbrv}

\end{document}